\documentclass{amsart}

\usepackage{amssymb}

\newtheorem{theorem}{Theorem}[section]
\newtheorem{lemma}[theorem]{Lemma}

\newtheorem{remark}[theorem]{Remark}

\numberwithin{equation}{section}
\newcommand{\Ai}{\text{Ai\,}}

\newcommand{\Tr}{\text{Tr\,}}
\newcommand{\re}{\text{Re\,}}

\begin{document}
\setcounter{page}{1}

\thanks{Supported by the G\"oran Gustafsson Foundation (KVA), the ESF-network MISGAM and
the EU network ENIGMA}

\dedicatory{To Percy Deift on his 60:th birthday}

\title[On some special directed last-passage percolation models]
{On some special directed last-passage percolation models}
\author[K.~Johansson]{Kurt Johansson}

\address{
Department of Mathematics,
Royal Institute of Technology,
SE-100 44 Stockholm, Sweden}

\email{kurtj@kth.se}

\begin{abstract}
We investigate extended processes given by last-passage times in directed models defined using
exponential variables with decaying mean. In certain cases we find the universal Airy process, but
other cases lead to non-universal and trivial extended processes.

\end{abstract}

\maketitle

\section{Introduction and results}\label{sect1}


Let $w(i,j)$ be independent exponential variables with parameter $t_i+t_j$, where $t_j>0$ are
given numbers, i.e.
\begin{equation}\label{1.0}
\mathbb{P}[w(i,j)\ge x]=e^{-(t_i+t_j)x},
\end{equation}
$i,j\ge 1$. 
In this paper we will consider the case when $t_i=i^\alpha$, $0<\alpha\le 1$.
Other models with varying parameters have been studied in \cite{GTW}.
Consider the last-passage times $G(m,n)$, $m,n\ge 1$, defined by
\begin{equation}\label{1.1}
G(m,n)=\max_{\pi}\sum_{(i,j)\in\pi} w(i,j),
\end{equation}
where the maximum is over all up/right paths $\pi$ from $(1,1)$ to $(m,n)$. 
This means that $\pi=((i_k,j_k))_{k=1}^{m+n-1}$, where $(i_1,j_1)=(1,1)$, $(i_{m+n-1},j_{m+n-1})
=(m,n)$ and $(i_{k+1},j_{k+1})-(i_k,j_k)=(0,1)$ or $(1,0)$.

It was proved in \cite{JoSh} that if $t_i+t_j=1$ for all $i,j\ge 1$, then for $m\ge n$,
\begin{equation}\label{1.2}
\mathbb{P}[G(m,n)\le\xi]=\frac 1{Z_{m,n}}\int_{[0,\xi]^n}\prod_{1\le i<j\le n}(x_i-x_j)^2\prod_{j=1}^n
x_j^{m-n}e^{-x_j}d^nx,
\end{equation}
which is the eigenvalue density in the Laguerre random matrix ensemble. This leads to the formula
\begin{equation}\label{1.2'}
\mathbb{P}[G(m,n)\le\xi]=\det(I-K_{m,n})_{L^2(\xi,\infty)}
\end{equation}
for the distribution function, \cite{Me}.
Here
\begin{equation}
K_{m,n}(x,y)=\sum_{j=0}^{n-1}p_j^{(m-n)}(x)p_j^{(m-n)}(y)(x^{m-n}e^{-x}y^{m-n}e^{-y})^{1/2},
\notag
\end{equation}
where $p_j^{(\alpha)}$, $j\ge 0$, are the normalized Laguerre polynomials.
The formula (\ref{1.2'}) was used in \cite{JoSh} to show that 
\begin{equation}\label{1.3}
\mathbb{P}[\frac{G(m,n)-cn}{dn^{1/3}}\le\xi]\to F_{\text{TW}}(\xi)
\end{equation}
as $n\to\infty$, $m/n\to\gamma\ge 1$, for some explicit constants $c,d$. Here $F_{\text{TW}}(\xi)=$
$\det(I-A)_{L^2(\xi,\infty)}$ is the Tracy-Widom distribution and
\begin{equation}
A(x,y)=\int_0^\infty\Ai(x+t)\Ai(y+t)dt
\notag
\end{equation}
is the Airy kernel.

In the general case with arbitrary $t_j>0$ there is no formula like (\ref{1.2}) which relates the
distribution of $G(m,n)$ to a random matrix ensemble. However we still have a formula like
(\ref{1.2'}) with an explicit correlation kernel, see (\ref{1.8}) below.

P. Forrester has noted that if we take $t_k=k+\beta$, $\beta>-1$, then we also have a limit law but 
with a different limiting distribution related to the distribution of the smallest eigenvalue in a 
Laguerre ensemble (hard edge limit):
\begin{equation}\label{1.4}
\mathbb{P}[G(n,n)-2\log n\le \xi]\to\det(I-K_\beta)_{L^2(\xi,\infty)}\doteq U_\beta(\xi)
\end{equation}
as $n\to\infty$. The kernel $K_\beta$ is related to the Bessel kernel,
\begin{equation}\label{1.5}
K_\alpha^{\text{Bessel}}(x,y)=\frac{J_\alpha(\sqrt{x})\sqrt{y}J'_\alpha(\sqrt{y})-
\sqrt{x}J'_\alpha(\sqrt{x})J_\alpha(\sqrt{y})}{2(x-y)},
\end{equation}
by
\begin{equation}\label{1.6}
\frac 1{\sqrt{uv}}K_\beta(\log\frac 1u,\log\frac 1v)=4K_{2\beta+1}^{\text{Bessel}}(4u,4v).
\end{equation}
Note that although the limit in (\ref{1.4}) is related to a universal distribution function from
random matrix theory the limit should be thought of as non-universal, since a small perturbation of 
the
distribution of the $w(i,j)$'s (perturbing $\beta$) changes the limiting distribution. In \cite{BoFo}
it is proved that
\begin{equation}
\lim_{\beta\to\infty} U_\beta(-2\log(4\beta)+(2\beta)^{-2/3}s)=F_{\text{TW}}(s).
\notag
\end{equation}
It can also be shown that, 
\begin{equation}
U_0(\xi)=e^{-e^{-\xi}}.
\notag
\end{equation}
This follows from formulas in \cite{FoRa}.
Hence the family $U_\beta(\xi)$ interpolates between the Gumbel and Tracy-Widom distributions,
compare \cite{JoGTW}.

We can also consider the process $k\to G(N+k,N-k)$, $|k|<N$. In the case when $t_i+t_j=1$ for all
$i,j\ge 1$ this process, appropriately rescaled in a neighbourhood of the origin, converges to 
the Airy process. In the case when the $w(i,j)$ are geometric random variables this is proved
in \cite{JoDPNG}, but the proof could be modified to the exponentual case. 
P. Forrester has raised the question \cite{Fo} what happens in the case $t_i=s_i=i+\beta\,$ ? 
What kind of extended limiting process do we get, an extended Bessel kernel process? In fact it 
turns out that we get a trivial extended process meaning that $G(N+k_1,N-k_1)$ and 
$G(N+k_2,N-k_2)$ have the same fluctuations for $k_1$ and $k_2$ far apart (of the order $N$).
One of the results of this paper is a proof of a weak version of this result. We will discuss the weak
version below, but first we will give a brief heuristic motivation why we can expect a trivial
extended process.

If we forget about the maximum in (\ref{1.1}) we are summing independent exponential random
variables with smaller and smaller variance, the parameter increasing linearly.
Let $X_j$, $j\ge 1$ be independent with distribution $\text{Exp}(j)$. Then
$$
\mathbb{E}[\sum_{j=m}^n X_j]=\sum_{j=m}^n\frac 1j\approx \log\frac nm,
$$
so we can expect a logarithm in the mean, which is exactly what we see in (\ref{1.4}). Consider
now the variance,
$$
\text{Var}[\sum_{j=m}^n X_j]=\sum_{j=m}^n\frac 1{j^2}\approx \frac 1m.
$$
This means that the contribution to the fluctuations should come from the $w(i,j)$'s with 
$i+j$ small. But these will contribute the same fluctuations to almost all the points on the line
$i+j=n$, when $n$ is large, and hence almost all points on this line should have the same fluctuations.

Note that if we instead take the $X_j$'s to be independent with distribution $\text{Exp}(j^\alpha)$,
$0<\alpha<1$, we get
\begin{equation}\label{1.7}
\text{Var}[\sum_{j=m}^n X_j]=\sum_{j=m}^n\frac 1{j^{2\alpha}}\approx 
\begin{cases} \frac 1{2\alpha-1}(m^{1-2\alpha}-n^{1-2\alpha}), &\frac 12<\alpha\le 1
\\ \frac 1{1-2\alpha}(n^{1-2\alpha}-m^{1-2\alpha}), &0<\alpha\le \frac 12,
\end{cases}
\end{equation}
so we can expect a difference between the cases $\alpha>1/2$ and $\alpha<1/2$.

The variance in (\ref{1.7}) really corresponds to moving along the axes, i.e. considering
$G(j,1)$ or $G(1,j)$, $j\ge 1$. We know from the case $t_i=1$ that as we move from the axes to the 
diagonal there is a reduction in the fluctuations exponent from $1/2$ to $1/3$, i.e. by $1/6$.
It is reasonable to expect a similar reduction in the caes when $t_i=i^\alpha$. Hence, the 
fluctuation exponent should be $\max(0,1/3-\alpha)$, which indicates a change when $\alpha =1/3$.
The above heuristics indicates that if we choose $t_i=i^\alpha$ we can expect changes in the behaviour
when $\alpha=1/2$ and $\alpha=1/3$. We will see below that this is indeed the case.

The probability measure which is the product measure of (\ref{1.0}) for $(i,j)\in\mathbb{Z}_+^2$,
$i+j\le 2N$, can be mapped to a determinantal point process on $\{-N+1,\dots,N-1\}\times\mathbb{R}$
with a last particle $x_{\max}^r$ on each line $\{r\}\times\mathbb{R}$, and such that 
$G(N+r,N-r)=x_{\max}^r$, $|r|<N$. This can be proved by modifying the argument for the geometric case
in sect. 5 of \cite{JoHou} or by taking an appropriate limit of the geometric case.
The resulting determinantal process has the correlation kernel
\begin{equation}\label{1.8}
K_N(r,x;s,y)=\tilde{K}_N(r,x;s,y)-\phi_{r,s}(x,y),
\end{equation}
where
\begin{equation}\label{1.9}
\tilde{K}_N(r,x;s,y)=\frac 1{(2\pi i)^2}\int_{\Gamma}dw\int_{\Gamma}dz e^{-xz-yw}\frac 1{z+w}F(z,w)
\end{equation}
and
\begin{equation}\label{1.10}
\phi_{r,s}(x,y)=\frac 1{2\pi}\int_{\mathbb{R}} e^{i\lambda(y-x)}F(e^{i\lambda},e^{i\lambda})d\lambda,
\end{equation}
if $r<s$ and $\phi_{r,s}\equiv 0$ if $r\ge s$. Here,
\begin{equation}\label{1.11}
F(z,w)=\frac{\prod_{k=1}^{N+r}(1+z/t_k)\prod_{k=1}^{N-s}(1+w/t_k)}
{\prod_{k=1}^{N-r}(1-z/t_k)\prod_{k=1}^{N+s}(1-w/t_k)},
\end{equation}
and $\Gamma=-\Gamma_- +\Gamma_+$, where $\Gamma_\pm$ are given by
$t\to te^{\pm i\pi/4}$, $t\ge 0$.

We want to consider scaling limits of the correlation kernel $K_N$ when $t_i=i+\beta$ and
$t_i=i^\alpha$, $0<\alpha<1$. Define
\begin{equation}\label{1.12}
c_{N,r}=\sum_{k=1}^{N+r}\frac 1{t_k}+\sum_{k=1}^{N-r}\frac 1{t_k},
\end{equation}
\begin{equation}\label{1.13}
G_{1,\beta}(z)=\prod_{k=1}^\infty\left(1-\frac z{k+\beta}\right)e^{z/k+\beta},
\end{equation}
\begin{equation}\label{1.14}
G_{\alpha}(z)=\prod_{k=1}^\infty\left(1-\frac z{k^\alpha}\right)e^{z/k^\alpha},
\end{equation}
for $1/2<\alpha<1$ and
\begin{equation}\label{1.15}
G_{\alpha}(z)=\prod_{k=1}^\infty\left(1-\frac z{k^\alpha}\right)e^{z/k^\alpha+z^2/2k^\alpha},
\end{equation}
for $1/3<\alpha\le 1/2$. We then have the following theorem.

\begin{theorem}\label{thm1.1}
a) If $t_i=i+\beta$, $\beta>-1$ or $t_i=i^\alpha$, $1/2<\alpha<1$, $i\ge 1$, then
\begin{equation}\label{1.16}
\tilde{K}_N(r,x+c_{N,r};s,y+c_{N,s})\to
\frac 1{(2\pi i)^2}\int_{\Gamma} dw\int_{\Gamma}dz\frac{e^{-xz-yw}}{z+w}\frac{G(-z)G(-w)}
{G(z)G(w)},
\end{equation}
uniformly for $x$ and $y$ in a compact set, as $N-|r|\to\infty$ and $N-|s|\to\infty$, where
$G=G_{1,\beta}$ and $G=G_\alpha$ respectively. Furthermore, for any $f$ such that $f$ and $\hat{f}$
belong to $L^1(\mathbb{R})$,
\begin{equation}\label{1.18}
\int_{\mathbb{R}} f(x)\phi_{r,s}(x+c_{N,r},y+c_{N,s})dx\to f(y)
\end{equation}
as $N-|r|\to\infty$ and $N-|s|\to\infty$.

b) Let $t_i=i^\alpha$, $i\ge 1$, $1/3<\alpha\le 1/2$. Set $r=[N\tanh\tau]$, $s=[N\tanh\sigma]$
if $\alpha=1/2$ and $r=[N\tau]$, $s=[N\sigma]$ if $1/3<\alpha<1/2$. Then,
\begin{equation}\label{1.19}
\tilde{K}_N(r,x+c_{N,r};s,y+c_{N,s})\to
\frac 1{(2\pi i)^2}\int_\Gamma dw\int_{\Gamma}dz\frac{e^{-xz-yw-\tau z^2+\sigma w^2}}{z+w}
\frac{G_\alpha(-z)G_\alpha(-w)}{G_\alpha(z)G_\alpha(w)},
\end{equation}
and
\begin{equation}\label{1.20}
\phi_{r,s}(x+c_{N,r},y+c_{N,s})\to\frac 1{\sqrt{4\pi (\sigma-\tau)}}e^{-(y-x)^2/4(\sigma-\tau)}
\end{equation}
uniformly for $x,y,\tau,\sigma$ in a compact set as $N\to\infty$.

c) Let $t_i=i^\alpha$, $0<\alpha\le 1/3$, $i\ge 1$. Define $d_N=(2\log N)^{1/3}$ if $\alpha=1/3$
and $d_N=2^{1/3}(1-3\alpha)^{-1/3}N^{1/3-\alpha}$, if $0<\alpha<1/3$. Then, with $r=[d_N^2N^{2\alpha}\tau]$,
$s=[d_N^2N^{2\alpha}\sigma]$,
\begin{equation}\label{1.21}
d_Ne^{\tau^3/3+\sigma^3/3-\sigma\eta+\tau\xi} K_N(r,[c_{N,r}+d_N(\xi-\tau^2)];r,
[c_{N,s}+d_N(\eta-\sigma^2)])\to A(\tau,\xi;\sigma,\eta)
\end{equation}
uniformly for $\sigma,\tau,\xi,\eta$ in a compact set as $N\to\infty$.
Here
\begin{equation}
A(\tau,\xi;\sigma,\eta)=\begin{cases}\int_0^\infty e^{-\lambda(\tau-\sigma)}\Ai(\xi+\lambda)
\Ai(\eta+\lambda)d\lambda, &\tau\ge\sigma \\
-\int_{-\infty}^0 e^{-\lambda(\tau-\sigma)}\Ai(\xi+\lambda)
\Ai(\eta+\lambda)d\lambda, &\tau<\sigma 
\end{cases}
\notag
\end{equation}
is the extended Airy kernel.
\end{theorem}

The proof of theorem \ref{thm1.1} will be given in section \ref{sect2.1}. It would be possible to 
consider more general sequences of parameters $(t_i)_{i\ge 1}$ with similar growth behaviour, but
we will not do that here.

\begin{remark} \rm Let $K_\beta(x,y)$ denote the right hand side of (\ref{1.16}), with $G=G_{1,\beta}$.
Set $\delta=2(\sum_{k=1}^\infty1/k(k+\beta)-\gamma)$, where $\gamma$ is Euler's constant. Then,
\begin{equation}\label{1.22}
\frac 1{\sqrt{uv}}K_\beta(\log\frac 1u+\delta,\log \frac 1v+\delta)=
4K^{\text{Bessel}}_{2\beta+1}(4u,4v), 
\end{equation}
which relates to the limit (\ref{1.4}).\it
\end{remark}

As we will argue below we should interpret the result a) in the theorem as saying that we have a
trivial extended process in the limit, i.e. corresponding points at different lines in the 
determinantal process will have identical fluctuations. Note that in case a) we have fluctuation exponent
0, transversal correlation exponent 1, and a non-universal limit in the sense that the limiting
correlation kernel depends on the $t_i$'s. In case b) we still have fluctuation exponent 0,
but the transversal correlation exponent varies as $2\alpha$, and we have a non-trivial and
non-universal limiting process.
Finally, in case c) we have a varying fluctuation exponent $1/3-\alpha$, the transversal correlation
exponent equals the standard value 2/3, and we have a universal limiting process, the Airy process.
In this case the specific choice of the $t_i$'s is not seen in the limit.

The limiting kernel given by (\ref{1.16}) and (\ref{1.18}) can be written as
\begin{equation}\label{1.23}
K(r,x;s,y)=K(x,y)-\delta(x,y)\eta_{rs},
\end{equation}
where $\eta_{rs}=1$ if $r<s$ and $=0$ otherwise.
Consider a point process that lives on two lines, denoted 1 and 2. Denote the particles on the first line
by $x_j$ and those on the second by $y_j$. That the processes on the two lines are identical 
should mean that
\begin{equation}\label{1.24}
\mathbb{E}[\prod_j(1-\phi_1(x_j))\prod_j(1-\phi_2(y_j))]=
\mathbb{E}[\prod_j(1-\phi_1(x_j)-\phi_2(x_j)+\phi_1(x_j)\phi_2(x_j))]
\end{equation}
for any continuous $\phi_1$, $\phi_2$, $0\le\phi_i\le 1$ with compact support. Let $z=(i,x)$, $i=1,2$,
be a point on any of the two lines, and set $\phi(z)=\phi(i,x)=\phi_i(x)$. Assume that 
the point process on both lines with points $\{z_i\}$ is determinatal with correlation kernel
\begin{equation}\label{1.25}
K_{\text{ext}}(i,x;j,y)=K(x,y)-\delta(x-y)\eta_{i,j}.
\end{equation}
Then (\ref{1.24}) will hold formally if
\begin{equation}\label{1.26}
\sum_{m=1}^\infty\frac{(-1)^{m-1}}{m}\Tr(K_{\text{ext}}\phi)^m=
\sum_{m=1}^\infty\frac{(-1)^{m-1}}{m}\Tr(K(\phi_1+\phi_2+\phi_1\phi_2))^m,
\end{equation}
since this is what we formally get if we take the logarithm of both sides in (\ref{1.24}) 
and use the basic properties of determinantal processes, \cite{Me}. 
The identity (\ref{1.26}) is, at least as
a formal identity, of combinatorial nature and will be discussed in section \ref{sect2.2}.

\section{Proofs}\label{sect2}

\subsection{Proof of theorem \ref{thm1.1}}\label{sect2.1}

We will use some notations, formulas and estimates from the theory of entire functions.
Set
$$
E(z;p)=(1-z)e^{z+z^2/2+\dots+z^p/p}
$$
and
$$
c_M^{(j)}=\sum_{k=1}^M\frac 1{t_k^j},
$$
so that $c_{N,r}=c_{N+r}^{(1)}+c_{N-r}^{(1)}$. Furthermore we write
$$
H_M(z)=\sum_{k=1}^M\log E(-\frac z{t_k};1)=\sum_{k=1}^M\log(1+z/t_k)-z/t_k.
$$
Introduce the counting function $n(t)=\#\{k\ge 1\,;\,t_k\le t\}$. We have that 
$n(t)=\max([t-\beta],0)$ if $t_i=i+\beta$ and $n(t)=[t^{1/\alpha}]$ if $t_i=i^\alpha$, 
where $[\cdot]$ denotes the integer part.
The following
estimate is standard
\begin{equation}\label{2.1}
|\log E(z;p)|\le\frac 1{1-q}|z|^{p+1}
\end{equation}
if $|z|\le q<1$, \cite{Bo}. Integration by parts gives the following identities
\begin{equation}\label{2.2}
H_M(z)=M\log E(-\frac z{t_M};1)-z^2\int_0^{t_M}\frac{n(t)}{t^2(z+t)}dt,
\end{equation}
\begin{equation}\label{2.3}
H_M(z)=M\log E(-\frac z{t_M};2)-\frac{c_M^{(2)}}2 z^2+ z^3\int_0^{t_M}\frac{n(t)}{t^3(z+t)}dt,
\end{equation}
\begin{equation}\label{2.4}
H_M(z)=M\log E(-\frac z{t_M};3)-\frac{c_M^{(2)}}2 z^2+\frac{c_M^{(3)}}3 z^3+ 
z^4\int_0^{t_M}\frac{n(t)}{t^4(z+t)}dt.
\end{equation}
It follows from (\ref{1.9}) that
\begin{align}\label{2.5}
&\tilde{K}_N(r,x+c_{N,r};s,y+c_{N,s})
\notag\\
&=\frac 1{(2\pi i)^2}\int_\Gamma dw\int_\Gamma dz\frac {e^{-xz-yw}}{z+w}e^{H_{N+r}(z)-H_{N-r}(-z)
+H_{N-s}(w)-H_{N+s}(-w)}.
\end{align}

Write $z=ue^{\pm i\pi/4}$, $u\ge 0$. Set $f(u)=\log(1+u^2+u\sqrt{2})-u\sqrt{2}$. 
Then, by (\ref{2.2}),
\begin{equation}\label{2.6}
\re H_M(ue^{\pm i\pi/4})=\frac M2 f(\frac u{t_M})-\frac {u^3}{\sqrt{2}}\int_0^{t_M}
\frac{n(t)}{t^2(t^2+u^2+ut\sqrt{2})}dt,
\end{equation}
for $u\ge 0$. Note that $f'(u)\le 0$ and $f(0)=0$, and hence $f(u)\le 0$ for $u\ge 0$ and $f(u)\ge 0$ for
$u\ge 0$. Together with (\ref{2.6}) this gives
\begin{align}\label{2.7}
&\re[-xz+H_{N+r}(z)-H_{N-r}(-z)]
\notag\\
&\le
-\frac{xu}{\sqrt{2}}-\frac{u^3}{\sqrt{2}}\left[\int_0^{t_{N+r}}\frac{n(t)}{t^2(t^2+u^2+ut\sqrt{2})}dt
+\int_0^{t_{N-r}}\frac{n(t)}{t^2(t^2+u^2-ut\sqrt{2})}dt\right].
\end{align}

\begin{lemma}\label{lem2.1}
Let $t_i=i+\beta$ or $t_i=i^\alpha$, $i\ge 1$, $0<\alpha<1$, $\beta>-1$. Set $\rho (u)=\log u$ and
$r(u)=u^{1/\alpha-1}$ respectively. There are constants $C$ and $D$, which only depend on $\alpha$
or $\beta$, so that, for all sufficiently large $M$,
\begin{equation}\label{2.8}
u^3\int_0^{t_M}\frac{n(t)}{t^2(t^2+u^2+ut\sqrt{2})}dt\ge Cu\rho (\min(t_M,u))
\end{equation}
if $u\ge D$.
\end{lemma}

\begin{proof}
Consider the case $t_i=i+\beta$, so that $n(t)=[t-\beta]$. Note that , if $t_M\ge 2(\beta+1)$, then
the left hand side of (\ref{2.8}) is
\begin{equation}\label{2.8.1}
\ge \frac{u^3}{2}\int_{2(\beta+1)}^{t_M}\frac{dt}{t(t^2+u^2+ut\sqrt{2})}
\end{equation}
If $u\ge t_M$, then $t\le t_M\le u$, and hence $t^2+u^2+ut\sqrt{2}\le(2+\sqrt{2})u^2$, and we 
see that the expression in (\ref{2.8.1}) is
$$
\ge\frac{u}{2(2+\sqrt{2})}\int_{2(\beta+1)}^{t_M}\frac{dt}{t}=\frac{u}{2(2+\sqrt{2})}
\log\frac{t_M}{2(\beta+1)}\ge Cu\rho (t_M)
$$
if $M$ is sufficiently large. In the case $t_i=i^\alpha$ we get similarly that the left hand side
of (\ref{2.8}) is
$$
\ge \frac{u}{2(2+\sqrt{2})}\int_{1}^{t_M}t^{1/\alpha-2}dt\ge Cu\rho (t_M)
$$
for all sufficientlty large $M$.
Assume now that $D\le u\le t_M$ with a suitable $D$. If $D\ge 2(\beta+1)$, then
the expression in (\ref{2.8.1}) is
$$
\ge\frac{u^3}{2}\int_{2(\beta+1)}^{u}\frac{dt}{t(t^2+u^2+ut\sqrt{2})}
\ge\frac{u}{2(2+\sqrt{2})}\int_{2(\beta+1)}^{u}\frac{dt}t\ge
Cu\log u,
$$
if we choose $D$ sufficiently large. The proof for $t_i=i^\alpha$ is completely analogous.
\end{proof}

Assume $t_i=i+\beta$ or $t_i=i^\alpha$ with $1/2<\alpha<1$, 
and that $x$ and $y$ belong to a compact set. It follows from (\ref{2.7})
that
\begin{equation}\label{2.9}
\log|e^{-xz+H_{N+r}(z)-H_{N-r}(-z)}|
\le -\frac{xu}{\sqrt{2}}-\sqrt{2} u^3\int_0^{t_{N-|r|}}\frac{n(t)}{t^2(t^2+u^2+ut\sqrt{2})}dt.
\end{equation}
Suppose that $|x/\sqrt{2}|\le K$ for some constant $K$. From lemma \ref{lem2.1} we see that 
the right hand side of (\ref{2.9}) is
\begin{equation}\label{2.9.2}
\le Ku-Cu\rho (\min(t_{N-|r|},u))
\end{equation}
if $u\ge D$ and $n-|r|$ is sufficiently large. We can thus choose $L$ so that if $u\ge L$
and $M$ is sufficiently large, then
\begin{equation}\label{2.10}
|e^{-xz+H_{N+r}(z)-H_{N-r}(-z)}|
\le e^{-u}
\end{equation}
for $z=ue^{\pm i\pi/4}$. The same type of estimate can be done for the $w$-part.

It follows from (\ref{2.1}), (\ref{2.2}) and $\alpha>1/2$ that
$$
\lim_{M\to\infty} e^{H_M(z)}=\prod_{k=1}^\infty(1+\frac z{t_k})e^{-z/t_k},
$$
uniformly on compact sets. From this and the estimate (\ref{2.10}) we now see that we can take the 
limit in the integral (\ref{2.5}) as $N-|r|$ and $N-|s|$ both tend to infinity and obtain the 
right hand side of (\ref{1.16}) with $G_{1,\beta}$ or $G_\alpha$.

We also want to prove (\ref{1.18}). Suppose $r<s$ and set
$$
\psi_{r,s}(t)=\frac 1{2\pi}\int_{\mathbb{R}} e^{i\lambda t}F_{N,r,s}(\lambda)d\lambda,
$$
where
\begin{equation}\label{2.10.2}
F_{N,r,s}(\lambda)=\prod_{k=N-s+1}^{N-r}\frac 1{E(i\lambda/t_k;1)}
\prod_{k=N+r+1}^{N+s}\frac 1{E(-i\lambda/t_k;1)},
\end{equation}
so that $\phi_{r,s}(x+c_{N,r},y+c_{N,s})=\psi_{r,s}(y-x)$. Thus,
\begin{equation}
\int_{\mathbb{R}} f(x)\phi_{r,s}(x+c_{N,r},y+c_{N,s})dx=
\frac 1{2\pi}\int_{\mathbb{R}}e^{i\lambda y}\hat{f}(\lambda)F_{N,r,s}(\lambda)d\lambda.
\notag
\end{equation}
Note that $|E(\pm i\lambda/t_k;1)|=(1+\lambda^2/t_k^2)^{1/2}\ge 1$ and consequently
$|F_{N,r,s}(\lambda)|\le 1$ for all $\lambda\in\mathbb{R}$. By dominated convergence it 
now suffices to show that $F_{N,r,s}(\lambda)\to 1$ pointwise as $N-|r|\to\infty$, 
$N-|s|\to\infty$, since we assume that $\hat{f}\in L^1(\mathbb{R})$.

It follows from (\ref{2.2}) that
\begin{align}\label{2.11}
&\sum_{k=N-s+1}^{N-r}\log E(i\lambda/t_k;1)=
\notag\\
&(N-r)\log E(i\lambda/t_{N-r};1)-
(N-s)\log E(i\lambda/t_{N-s};1)+\lambda^2\int_{t_{N-s}}^{t_{N-r}}\frac{n(t)}{t^2(t-i\lambda)}dt.
\end{align}
For a fixed $\lambda$ and $N-r$ large it follows from (\ref{2.1}) that
$$
(N-r)|\log E(i\lambda/t_{N-r};1)|\le 2\lambda^2\frac{N-r}{t_{N-r}^2}=2\lambda^2 (N-r)^{1-2\alpha},
$$
which $\to 0$ as $N-s\to \infty$ since $\alpha>1/2$. The second term in (\ref{2.11}) is treated similarly.
The third term is estimated as follows
$$
\left|\int_{t_{N-s}}^{t_{N-r}}\frac{n(t)}{t^2(t-i\lambda)}dt\right|\le
\int_{t_{N-s}}^{t_{N-r}}\frac{n(t)}{t^3}dt\to 0
$$
as $N-s\to\infty$, since $n(t)=[t-\beta]$ or $n(t)=[t^{1/\alpha}]$ with $\alpha>1/2$. Similarly, we
can show that
$$
\sum_{k=N+r+1}^{N+s}\log E(-i\lambda/t_k;1)\to 0
$$
as $N-|r|\to\infty$, $N-|s|\to\infty$.

Consider next the case when $1/3<\alpha\le 1/2$.

Assume that $x$ and $y$ belong to a compact set and that $r=[\tau N^{2\alpha}]$, 
$s=[\sigma N^{2\alpha}]$, 
or $r=[N\tanh\tau]$, $s=[N\tanh \sigma]$ in case $\alpha=1/2$, 
where $\sigma, \tau$ belong to a compact set. We use again
the formula (\ref{2.5}) for $\tilde{K}_N$ and the estimate (\ref{2.9}). We see that $N-|r|\ge N/2$ if
$N$ is sufficiently large in case $1/3<\alpha\le 1/2$. the expression in (\ref{2.9}) is 
bounded by (\ref{2.9.2}) by lemma \ref{lem2.1}. with $r(u)=u^{1/\alpha-1}$, and again we obtain 
the estimate (\ref{2.10}) for $z=u e^{\pm i\pi/4}$ and $u\ge D$ with a suitable $D$. We can write
$$
E^{H_M(z)}=\prod_{k=1}^M E(-z/t_k;2) e^{-c_M^{(2)}z^2}.
$$
Hence,
$$
e^{H_{N+r}(z)-H_{N-r}(-z)}=\frac{\prod_{k=1}^{N+r} E(-z/t_k;2)}{\prod_{k=1}^{N-r} E(z/t_k;2)}
e^{(c_{N-r}^{(2)}-c_{N+r}^{(2)})z^2}.
$$
Since $\sum_{k=1}^\infty 1/t_k^3<\infty$ it follows from the theory of canonical products that
$$
\lim_{M\to\infty}\prod_{k=1}^ME(\pm z/t_k;2)=\prod_{k=1}^{\infty}E(\pm z/t_k;2)=G_\alpha(\mp z),
$$
uniformly on compacts. Now, if $1/3<\alpha<1/2$,
\begin{align}
c_{N-r}^{(2)}-c_{N+r}^{(2)}&=-\text{sgn\,}(r)\sum_{k=N-|r|+1}^{N+|r|}\frac 1{k^{2\alpha}}
\notag\\
&\sim -\text{sgn\,}(r)\frac{N^{1-2\alpha}}{1-2\alpha}\left[(1+|r|/N)^{1-2\alpha}-
(1-|r|/N)^{1-2\alpha}\right]\to -2\tau
\notag
\end{align}
as $N\to\infty$. The case $\alpha=1/2$ is analogous. This proves (\ref{1.19}). 

We also want to show (\ref{1.20}). Suppose that $r<s$ and
consider $\phi_{r,s}$. We can write
$$
\phi_{r,s}(x+c_{N,r},y+c_{N,s})=\psi_{r,s}(y-x),
$$
where
$$
\psi_{r,s}(t)=\frac 1{2\pi}\int_{\mathbb{R}} e^{i\lambda t}F_{N,r,s}(\lambda)d\lambda,
$$
and
\begin{align}
&F_{N,r,s}(\lambda)
\notag\\
&=\prod_{k=N-s+1}^{N-r}\frac 1{E(i\lambda/t_k;2)}
\prod_{k=N+r+1}^{N+s}\frac 1{E(-i\lambda/t_k;2)}
\prod_{k=N-s+1}^{N-r}e^{-\lambda^2/t_k^2}
\prod_{k=N+r+1}^{N+s}e^{-\lambda^2/t_k^2}.
\notag
\end{align}
From the convergence of the canonical products we see that 
$$
\lim_{N\to\infty}\prod_{k=N-s+1}^{N-r}\frac 1{E(i\lambda/t_k;2)}
\prod_{k=N+r+1}^{N+s}\frac 1{E(-i\lambda/t_k;2)}=1
$$
for each $\lambda\in\mathbb{R}$. If $1/3<\alpha<1/2$, then
$$
\sum_{k=N-s+1}^{N-r}\frac 1{k^{2\alpha}}+\sum_{k=N+r+1}^{N+s}\frac 1{k^{2\alpha}}
\sim\frac{2(s-r)}{N^{2\alpha}}\sim 2(\sigma-\tau)
$$
as $N\to\infty$. Since $\sigma>\tau$ it follows that
$$
\lim_{N\to\infty} \psi_{[\tau N^{2\alpha}],[\psi N^{2\alpha}]}(t)=\frac 1{2\pi}\int_{\mathbb{R}}
e^{i\lambda t-(\sigma-\tau)\lambda^2}d\lambda
=\frac 1{\sqrt{4\pi (\sigma-\tau)}} e^{-t^2/4(\sigma-\tau)},
$$
and we have proved (\ref{1.20}). In the case $\alpha=1/2$ we similarly get (\ref{1.20})
using the new expressions of $r,s$ in terms of $\tau,\sigma$.

It remains to treat the caes $0<\alpha<1/3$. Again our starting point is the formula (\ref{2.5})
and we will use the estimate (\ref{2.9}).
However, we need a new estimate of the integral in (\ref{2.9}).

\begin{lemma}\label{lem2.2}
Assume $0<\alpha\le 1/3$. If $u\ge t_M$, then
\begin{equation}\label{2.12}
u^3\int_0^{t_M}\frac{n(t)}{t^2(t^2+u^2+ut\sqrt{2})} dt\ge Cu(M^{1-\alpha}-1)
\end{equation}
for some constant $C>0$ that only depends on $\alpha$. If $0<\alpha<1/3$, there is a constant $C$, 
which only depends on $\alpha$, such that for $0\le u\le t_M$,
\begin{equation}\label{2.13}
u^3\int_0^{t_M}\frac{n(t)}{t^2(t^2+u^2+ut\sqrt{2})} dt\ge C\left[\frac{u^{1/\alpha}-u}{1-\alpha}
+u^3\frac{M^{1-3\alpha}-\max(u,1)^{1/\alpha-3}}{1-3\alpha}\right].
\end{equation}
If $\alpha=1/3$, there is a constant $C$ such that for $0\le u\le t_M$,
\begin{equation}\label{2.14}
u^3\int_0^{t_M}\frac{n(t)}{t^2(t^2+u^2+ut\sqrt{2})} dt\ge Cu^3\log\frac{t_M}{\max(u,1)}.
\end{equation}
\end{lemma}

\begin{proof}
If $u\ge t_M$, then using $t^2+u^2+ut\sqrt{2}\le (2+\sqrt{2})u^2$ we get
$$
u^3\int_0^{t_M}\frac{n(t)}{t^2(t^2+u^2+ut\sqrt{2})} dt\ge 
\frac{u}{1+\sqrt{2}}\int_1^{t_M}\frac{[t^{1/\alpha}]}{t^2}dt,
$$
which gives (\ref{2.12}).
If $0\le u\le t_M$, we write the left hand sides of (\ref{2.13}) and (\ref{2.14}) as
$$
u^3\int_0^{u}\frac{n(t)}{t^2(t^2+u^2+ut\sqrt{2})} dt
+u^3\int_u^{t_M}\frac{n(t)}{t^2(t^2+u^2+ut\sqrt{2})} dt.
$$
In the first integral we use again $t^2+u^2+ut\sqrt{2}\le (2+\sqrt{2})u^2$ and $n(t)=0$ for
$0\le t\le 1$, and in the second integral we use
$t^2+u^2+ut\sqrt{2}\le (2+\sqrt{2})t^2$. This yields immedeiately the estimates in the lemma.
\end{proof}

We now consider the case $0<\alpha<1/3$.
Write $c_0=2^{1/3}(1-3\alpha)^{-1/3}$ and let $r=[c_0^2\tau N^{2/3}]$, $s=[c_0^2\sigma N^{2/3}]$,
$x=[c_0N^{1/3-\alpha}(\xi-\tau^2)]$ and $y=[c_0N^{1/3-\alpha}(\eta-\tau^2)]$. 
Below we will ignore the fact that we take the integer parts.
In (\ref{2.5}) we
do the rescaling $z=c_0^{-1} N^{\alpha-1/3}\zeta$, $w=c_0^{-1} N^{\alpha-1/3}\omega$.
We will write $\xi'=\xi-\tau^2$, $\eta'=\eta-\sigma^2$. The integral in (\ref{2.5}) becomes
\begin{equation}\label{2.14.2}
\frac{c_0^{-1}N^{\alpha-1/3}}{(2\pi i)^2}\int_{\Gamma}d\omega\int_{\Gamma}d\zeta
\frac{e^{-\xi'\zeta-\eta'\omega}}{\zeta+\omega}e^{H_{N+r}(z)-H_{N-r}(-z)+H_{N-s}(w)-H_{N+s}(-w)}.
\end{equation}
The estimate (2.9) becomes, $z=ue^{\pm i\pi/4}\doteq c_0^{-1}N^{\alpha-1/3} ve^{\pm i\pi/4}$,
\begin{align}\label{2.15}
\log\left|e^{-\xi'\zeta+H_{N+r}(z)-H_{N-r}(-z)}\right|&\le 
-\frac{\xi' v}{\sqrt{2}}-\sqrt{2} u^3\int_0^{t_{N-|r|}}\frac{n(t)}{t^2(t^2+u^2+ut\sqrt{2})}dt
\notag\\
&\le 
-\frac{\xi' v}{\sqrt{2}}-\sqrt{2} u^3\int_0^{(N/2)^\alpha}\frac{n(t)}{t^2(t^2+u^2+ut\sqrt{2})}dt
\end{align}
if $N$ is sufficiently large. If $u\ge (N/2)^{\alpha}$, then the last expression in (\ref{2.15}) is
$$
\le -\xi' v/\sqrt{2}v-Cc_0^{-1}N^{\alpha-1/3}v((N/2)^{1-\alpha}-1)\le -C'N^{2/3} v
$$
for some constant $C'>0$ if $N$ is sufficiently large and $\xi,\tau$ belong to a compact set.
If $(N/4)^{\alpha}\le u\le (N/2)^{\alpha}$, then the last expression in (2.15) is
$$
\le  -\xi' v/\sqrt{2}-C(1-\alpha)^{-1} u(u^{1/\alpha}-1)\le -C'N^{2/3} v,
$$
and if $0\le u\le (N/4)^{\alpha}$ it follows from (\ref{2.13}) that the last expression in 
(\ref{2.15}) is $\le -\xi' v/\sqrt{2}v-C'N^{1-3\alpha}u^3\le -\xi' v/\sqrt{2}v-C''v^3$.
It follows from these estimates, that we can restrict the integration in (\ref{2.14.2}) to
$0\le v\le N^\gamma$, for any $\gamma>0$, with a negligible error in the limit $N\to\infty$.

By (\ref{2.14}) we can write
\begin{align}\label{2.15.2}
H_{N+r}(z)-H_{N-r}(-z)&=(N+r)\log E(-z/t_{N+r};3)-(N-r)\log E(z/t_{N-r};3)
\notag\\
&-\frac 12(c_{N+r}^{(2)}-c_{N-r}^{(2)})z^2+\frac 13(c_{N+r}^{(3)}-c_{N-r}^{(3)})z^3
\notag\\
&+ z^4\int_0^{t_{N+r}}\frac{n(t)}{t^4(t+z)}dt-
z^4\int_0^{t_{N-r}}\frac{n(t)}{t^4(t-z)}dt.
\end{align}
We choose $0<\gamma<\min(1/12,1/3-\alpha)$. Now, $c_{N+r}^{(2)}-c_{N-r}^{(2)}\sim 2r/N^{2\alpha}$ and
$c_{N+r}^{(3)}-c_{N-r}^{(3)}\sim c_0^3N^{1-3\alpha}$ as $N\to\infty$.
Hence,
$$
-\frac 12(c_{N+r}^{(2)}-c_{N-r}^{(2)})z^2+\frac 13(c_{N+r}^{(3)}-c_{N-r}^{(3)})z^3
\to-\tau\zeta^2+\frac 13\zeta^3
$$
as $N\to\infty$. Note that $|z/t_{N+r}|\le CN^{\gamma-1/3}<1/2$ if $N$ is large enough. Hence, it follows
from (\ref{2.1}) that 
$$
|(N+r)\log E(-z/t_{N+r};3)|\le 2(N+r)|-z/t_{N+r}|^4\le CN^{4\gamma-1/3},
$$
which $\to 0$ as $N\to\infty$ since $\gamma<1/12$.

We have, for $t\ge 1$ and $N$ sufficiently large,
$$
|t+z|\ge |t|-|z|\ge \frac 12|t|+\frac 12-c_0^{-1}N^{\alpha-1/3}|v|\ge \frac 12 |t|,
$$
since $N^{\alpha-1/3}|v|\le N^{\alpha-1/3+\gamma}\to 0$ as $N\to\infty$. 
Because $n(t)=0$ if $t<1$, we see that
$$
\left|z^4\int_0^{t_{N+r}}\frac{n(t)}{t^4(t+z)}dt\right|\le CN^{4\alpha-4/3+4\gamma}
\int_1^{(N+r)^\alpha}t^{1/\alpha-5}dt,
$$ 
which $\to 0$ as $N\to\infty$. The second integral in (\ref{2.15.2}) is analogous, and we
can treat $H_{N-s}(w)-H_{N+s}(-w)$ in exactly the same way.

We have proved that 
\begin{align}\label{2.16}
&\lim_{N\to\infty}c_0 N^{1/3-\alpha}\tilde{K}_N([c_0^2\tau N^{2/3}],[c_{N,r}+c_0N^{1/3-\alpha}
(\xi-\tau^2)];
\notag\\
&[c_0^2\sigma N^{2/3}],[c_{N,r}+c_0N^{1/3-\alpha}(\eta-\sigma^2)])
\notag\\
&=\frac 1{(2\pi i)^2}\int_\Gamma d\omega\int_\Gamma d\zeta\frac{e^{-(\xi-\tau^2)\zeta-
(\eta-\sigma^2)\omega}}{\zeta+\omega} e^{-\tau\zeta^2+\zeta^3/3+\sigma\omega^2+\omega^3/3}.
\end{align}

Note now that we have the identity
\begin{align}\label{2.17}
&\frac 1{(2\pi i)^2}\int_\Gamma d\omega\int_\Gamma d\zeta\frac{e^{-\xi\zeta-
\eta\omega}}{\zeta+\omega} e^{-\tau\zeta^2+\zeta^3/3+\sigma\omega^2+\omega^3/3}
\notag\\
&= e^{2(\sigma^3-\tau^3)/3+\sigma\eta-\xi\tau}\int_0^\infty e^{(\sigma-\tau)\lambda}
\Ai (\xi+\tau^2+\lambda)\Ai(\eta+\sigma^2+\lambda)d\lambda.
\end{align}
To see this observe that for $\re (\zeta+\omega)>0$ we have
$$
\int_0^\infty e^{-\lambda(\zeta+\omega)}d\lambda=\frac 1{\zeta+\omega}.
$$
Hence, the left hand side of (\ref{2.17}) can be written
$$
\int_0^\infty\left(\frac 1{2\pi i}\int_\Gamma e^{-\zeta(\xi+\lambda)-\tau\zeta^2+\zeta^3/3}d\zeta
\right)\left(\frac 1{2\pi i}\int_\Gamma e^{-\omega(\eta+\lambda)-\sigma\omega^2+\omega^3/3}d\omega
\right)d\lambda.
$$
Let $\Gamma'$ consist of the two rays $-\Gamma_1'$ and $\Gamma_2'$, where $\Gamma_1':t\to t
e^{3\pi i/4}$ and $\Gamma_1':t\to t e^{\pi i/4}$, $t\ge 0$. The change of variables
$\zeta=-iz$ gives
\begin{align}
\frac 1{2\pi i}\int_\Gamma e^{-\zeta(\xi+\lambda)-\tau\zeta^2+\zeta^3/3}d\zeta
&=\frac 1{2\pi i}\int_{\Gamma'} e^{iz(\xi+\lambda)+\tau z^2+z^3/3}dz
\notag\\
&=e^{-2\tau^3/3-(\xi+\lambda)\tau}\Ai(\xi+\lambda+\tau^2),
\notag
\end{align}
and we obtain (\ref{2.17}).

Suppose next that $r<s$ and consider $\phi_{r,s}$. Set $d_N=c_0N^{1/3-\alpha}$. Then
$$
d_N\psi_{r,s}(d_Nt)=\frac{1}{2\pi}\int_{\mathbb{R}}e^{i\lambda t}F_{N,r,s}(\frac{\lambda}{d_N})
d\lambda,
$$
where $F_{N,r,s}$ is given by (\ref{2.10.2}).
We have that
$$
\log F_{N,r,s}(\frac{\lambda}{d_N})=-H_{N-r}(-\frac{i\lambda}{d_N})+H_{N-s}(-\frac{i\lambda}{d_N})
-H_{N+s}(\frac{i\lambda}{d_N})+H_{N+r}(\frac{i\lambda}{d_N}).
$$
It follows from (\ref{2.15.2}) and appropriate estimates similar to the ones above that, for
$\sigma>\tau$,
$$
\lim_{N\to\infty}F_{N,r,s}(\frac{\lambda}{d_N})=-(\sigma-\tau)\lambda^2
$$
and also that we have an estimate
$$
\left|F_{N,r,s}(\frac{\lambda}{d_N})\right|\le Ce^{(\sigma-\tau)\lambda^2/2}
$$
if $N$ is sufficiently large. Thus
$$
\lim_{N\to\infty} d_N\psi_{r,s}(d_Nt)=\frac 1{\sqrt{4\pi(\sigma-\tau)}}e^{-t^2/4(\sigma-\tau)}.
$$
Since, $\phi_{r,s}(x+c_{N,r},y+c_{N,s})=\psi_{r,s}(y-x)$, we obtain
\begin{align}\label{2.18}
&\lim_{N\to\infty}\phi_{[c_0^2\tau N^{2/3}],[c_0^2\tau N^{2/3}]}([c_N,r+c_0 N^{1/3-\alpha}\xi],
[c_N,r+c_0 N^{1/3-\alpha}\eta])
\notag\\
&=\frac 1{\sqrt{4\pi(\sigma-\tau)}}e^{-t^2/4(\sigma-\tau)}.
\end{align}
If we combine (\ref{2.16})-(\ref{2.18}) we get (\ref{1.21}), since, \cite{Ok},
$$
\int_{\mathbb{R}}e^{-\lambda(\tau-\sigma)}\Ai(\xi+\lambda)\Ai(\eta+\lambda)d\lambda
=\frac 1{\sqrt{4\pi(\sigma-\tau)}}e^{-(\xi-\eta)^2/4(\sigma-\tau)-(\sigma-\tau)(\xi+\eta)
(\sigma-\tau)^3/12}.
$$

The case $\alpha=1/3$ is treated similarly. We replace $d_N=c_0 N^{1/3-\alpha}$ with $(2\log N)^{1/3}$ 
and $N^\gamma$ with $(\log N)^{1/4}$.

\subsection{The identity (\ref{1.26})}\label{sect2.2}

In this section we will discuss the identity (\ref{1.26}).
Write $g=\phi_1+\phi_2+\phi_1\phi_2$. We have
\begin{equation}\label{2.19}
\Tr (Kg)^m=\int_{\mathbb{R}^m} g(x_1)K(x_1,x_2)g(x_2)K(x_2,x_3)\dots g(x_m)K(x_m,x_1)d^mx
\end{equation}
and
\begin{align}\label{2.20}
\Tr (K_{\text{ext}}\phi)^m=\sum_{i_1,\dots,i_m=1,2}&\int_{\mathbb{R}^m} 
\phi_{i_1}(x_1)K_{i_1,i_2}(x_1,x_2)\phi_{i_2}(x_2)K_{i_2,i_3}(x_2,x_3)
\notag\\
&\dots \phi_{i_m}(x_m)
K_{i_m,i_1}(x_m,x_1)d^mx,
\end{align}
where we have written
\begin{equation}\label{2.21}
K_{ij}(x,y)=K(x,y)-\delta(x-y)\eta_{ij}.
\end{equation}
If we insert $g=\phi_1+\phi_2+\phi_1\phi_2$ into (\ref{2.19}) we can at each position 
choose $\phi_1,\phi_2$ or $\phi_1\phi_2$. Choosing $\phi_1$ or $\phi_2$ corresponds exactly to 
the summation over $i=1,2$ for a factor $\phi_i$ in (\ref{2.20}). If we choose $\phi_1\phi_2$, this must
correspond to the $\delta$-function contribution when we insert (\ref{2.21}) into (\ref{2.20}).
Note that when we insert (\ref{2.21}) into the product in (\ref{2.20}) we do not get any 
contribution from two consecutive $\delta$-functions, since $\eta_{i_1i_2}\eta_{i_2i_3}=0$ for all
choices of $i_1,i_2,i_3$. This means that we only get pairs $\phi_1\phi_2$ in (\ref{2.20})
just as in (\ref{2.19}). If we have $r$ factors of the type $\phi_1\phi_2$ in (\ref{2.19}) 
we can place them in the product in $\binom{m}{r}$ different ways. This must come from taking $r$
$\delta$-functions in $\Tr(K_{\text{ext}}\phi)^{m+r}$. 

Hence to show that the left and right hand sides of (\ref{1.26}) are equal we must show
that the signs and the combinatorial factors agree. The signs are easy, since the sign 
in front of $\Tr(K_{\text{ext}}\phi)^{m+r}$ is $(-1)^{m+r-1}$ and the sign from the
$\delta$-functions is $(-1)^r$. This gives $(-1)^{m-1}$ which is the sign in front of
$\Tr(Kg)^m$.

The coefficient in front of the expansion of the expression (\ref{2.19}) with $r$ factors 
$\phi_1\phi_2$ is $\frac 1m\binom {m}{r}$. The factor in front of $\Tr(K_{\text{ext}}\phi)^{m+r}$
is $1/(m+r)$ and hence we must show that the number of ways of placing the $r$ $\delta$-functions
must be
\begin{equation}\label{2.22}
\frac{m+r}m \binom{m}{r}=\binom{m}{r}+\binom{m-1}{r-1}.
\end{equation}
Recall from above that we cannot have two consecutive $\delta$-functions in the expression 
(\ref{2.20}) for the trace since this gives a zero contribution. The circular structure of 
the trace means that we have the following combinatorial problem: Coose $r$ points on the discrete
circle with $m+r$ points in such a way that the distances between the chosen points are all $\ge 2$.
We must show that the number of ways this can be done equals the expression in (\ref{2.22}).

Number the points as $0,1,\dots,m+r-1$ and count modula $m+r$. Let $c$ be the first point that
is included and let $\ell_1,\dots,\ell_{r-1}$ be the distances between the included points.
We get two contributions.

1) If $c=0$ we get
$$
\sum_{\ell_1+\dots+\ell_{r-1}\le m+r-2, \ell_i\ge 2} 1,
$$
since we cannot choose the last point.

2) If $c\neq 0$ we get
$$
\sum_{c=1}^{m-r+1}\sum_{\ell_1+\dots+\ell_{r-1}\le m+r-1-c, \ell_i\ge 2} 1.
$$
Write $\ell_i=k_i+2$. Then
\begin{equation}\label{2.23}
\sum_{\ell_1+\dots+\ell_{r-1}\le m+r-1-c, \ell_i\ge 2} 1
=\sum_{k_1+\dots+k_{r-1}\le m-r+1-c,k_i\ge 0} 1.
\end{equation}
We now use the identity
$$
\sum_{k_1+\dots+k_{r-1}=p,k_i\ge 0} 1=\binom{r+p-2}{r-2}.
$$
Hence, the expression in (\ref{2.23}) equals
$$
\sum_{p=0}^{m-r+1-c}\binom{r+p-2}{r-2}=\binom{m-c}{r-1}
$$
by the identity
\begin{equation}\label{2.24}
\sum_{k=0}^m\binom{n+k}{n}=\binom{n+m+1}{n+1}.
\end{equation}
This now shows that the expression in 1) is $\binom{m-1}{r-1}$.
The expression in 2) equals
$$
\sum_{c=1}^{m-r+1}\binom{m-c}{r-1}=\binom{m}{r}
$$
by (\ref{2.24}). This completes the argument.

\medskip
\noindent
{\bf Acknowledgement}: I thank Peter Forrester for bringing up the problem on the
properties of extended processes when we have decaying parameters during a visit
to Stockholm and for interesting discussions.

\end{document}